\documentclass[preprint,11pt]{elsarticle}

\usepackage{amssymb}
\usepackage{latexsym}
\usepackage{amsmath}
\usepackage{amsfonts}
\usepackage{amssymb}
\usepackage{amsthm}
\usepackage{psfrag}
\usepackage{epsfig}

\newsavebox{\fmbox}

\newtheorem{theorem}{Theorem}[section]
\newtheorem{lemma}[theorem]{Lemma}
\newtheorem{corollary}[theorem]{Corollary}
\newtheorem{remark}[theorem]{Remark}
\newtheorem{definition}[theorem]{Definition}
\newtheorem{proposition}[theorem]{Proposition}
\newtheorem{example}[theorem]{Example}

\newcommand{\eps}{\varepsilon}

\def\Vo{\vbox{\offinterlineskip\hbox{\kern 3pt$\scriptstyle\circ$}
\kern 1pt\hbox{$V$}}}
\def\Ho{\vbox{\offinterlineskip\hbox{\kern 3pt$\scriptstyle\circ$}
\kern 1pt\hbox{$H$}}}
\def\Wo{\vbox{\offinterlineskip\hbox{\kern 3pt$\scriptstyle\circ$}
\kern 1pt\hbox{$W$}}}

\let\vv =\v

\newcommand{\be}{\begin{equation}}
\newcommand{\ee}{\end{equation}}
\newcommand{\beq}{\begin{eqnarray}}
\newcommand{\eeq}{\end{eqnarray}}
\newcommand{\beqs}{\begin{eqnarray*}}
\newcommand{\eeqs}{\end{eqnarray*}}
\newcommand{\bt}{\begin{theorem}}
\newcommand{\et}{\end{theorem}}
\newcommand{\bex}{\begin{example}}
\newcommand{\eex}{\end{example}}
\newcommand{\br}{\begin{remark}}
\newcommand{\er}{\end{remark}}
\newcommand{\bc}{\begin{corollary}}
\newcommand{\ec}{\end{corollary}}
\newcommand{\bl}{\begin{lemma}}
\newcommand{\el}{\end{lemma}}
\newcommand{\bp}{\begin{proposition}}
\newcommand{\ep}{\end{proposition}}
\newcommand{\bd}{\begin{definition}}
\newcommand{\ed}{\end{definition}}

\newcommand{\R}{\mathbb{R}}

\newcommand{\PP}{\mathbb{P}}

\newcommand{\caI}{{\cal I}}

\newcommand{\caR}{{\cal R}}

\newcommand{\caP}{{\cal P}}

\newcommand{\caC}{{\cal C}}

\def\u{{\bf u}}

\def\n{{\bf n}}

\def\f{{\bf f}}
\def\g{{\bf g}}

\def\vv{{\bf v}}
\def\w{{\bf w}}

\def\x{{\bf x}}

\def\V{{\bf V}}
\def\K{{\bf K}}

\def\0{{\bf 0}}

\newcommand{\bnabla}{\mbox{\boldmath{$\nabla$}}}

\newcommand{\bs}{\mbox{\boldmath{$\sigma$}}}
\newcommand{\e}{\mbox{\boldmath{$\varepsilon$}}}

\newcommand{\Gc}{\Gamma_{C}}

\newcommand{\dg}{d\Gamma}

\newcommand{\restrictiona}[1]
{{ \begin{picture}(13,10) \put(-1,-4){$\mid$} \end{picture}
\hspace{-1em}_{_{\mbox{$_{#1}$}}}
}}
\newcommand{\cqfd}{{$\mbox{}$\hfill$\square$}}

\newcommand{\BB}{\mathbf{B}}

\journal{Applied Numerical Mathematics}

\begin{document}

\begin{frontmatter}

%% Title, authors and addresses

%% use the tnoteref command within \title for footnotes;
%% use the tnotetext command for the associated footnote;
%% use the fnref command within \author or \address for footnotes;
%% use the fntext command for the associated footnote;
%% use the corref command within \author for corresponding author footnotes;
%% use the cortext command for the associated footnote;
%% use the ead command for the email address,
%% and the form \ead[url] for the home page:
%%
%% \title{Title\tnoteref{label1}}
%% \tnotetext[label1]{}
%% \author{Name\corref{cor1}\fnref{label2}}
%% \ead{email address}
%% \ead[url]{home page}
%% \fntext[label2]{}
%% \cortext[cor1]{}
%% \address{Address\fnref{label3}}
%% \fntext[label3]{}

\title{On convergence of the penalty method for unilateral contact problems}

%% use optional labels to link authors explicitly to addresses:
%% \author[label1,label2]{<author name>}
%% \address[label1]{<address>}
%% \address[label2]{<address>}

\author[fc]{Franz Chouly\corref{fcc}}
\ead{franz.chouly@univ-fcomte.fr}

\author[ph]{Patrick Hild}
\ead{patrick.hild@math.univ-toulouse.fr}

\cortext[fcc]{Corresponding author. Phone: +33 3 81 66 64 89. Fax : +33 3 81 66 66 23.}

\address[fc]
{Laboratoire de Math\'ematiques - UMR CNRS 6623,
Universit\'e de Franche Comt\'e, 16 route de Gray, 25030 Besan\c{c}on Cedex,
France.
%
%email: {\tt franz.chouly@univ-fcomte.fr}
}

\address[ph]{
Institut de Math\'ematiques de Toulouse - UMR 5219 (CNRS/INSAT/UT1/UT2/UT3),
Universit\'e Paul Sabatier (UT3), 118 route de Narbonne, 31062 Toulouse Cedex 9,
France.
%
%email: {\tt patrick.hild@math.univ-toulouse.fr}
}

\begin{abstract}
%% Text of abstract
We present a convergence analysis of the penalty method applied to unilateral contact problems in two and three space dimensions. We first consider, under various regularity assumptions on the exact solution to the unilateral contact problem, the convergence of the continuous penalty solution as the penalty parameter $\varepsilon$ vanishes. Then, the analysis of the finite element discretized penalty method is carried out. Denoting by $h$ the discretization parameter, we show that the error terms we consider give the same estimates as in the case of the constrained problem when the penalty parameter is such that $\varepsilon = h$.
\end{abstract}

\begin{keyword}
%% keywords here, in the form: keyword \sep keyword
unilateral contact, variational inequality, finite elements,
penalty method, a priori error estimates.

%% MSC codes here, in the form: \MSC code \sep code
%% or \MSC[2008] code \sep code (2000 is the default)
\end{keyword}

\end{frontmatter}

%%
%% Start line numbering here if you want
%%
% \linenumbers

%% main text

%%%%%%%%%%%%%%%%%%%%%%%%%%%%%%%%%%%%%%%%%%%%%%%%
\section{Introduction}
\label{sect:introduction}

The penalty method is a classical and widespread method for the numerical treatment of constrained problems, in particular the unilateral contact problems arising in mechanics of deformable bodies which involve a nonlinear boundary condition written as an inequality (see, e.g.,  \cite{kikuchi-oden-88, laursen-02, wriggers-02}). Nevertheless, and to the best of our knowledge, the convergence analysis of the method in the simplest case of linear elastostatics with or without finite element discretization has been object of few studies. We may nevertheless quote the earlier, and pioneering works of Kikuchi,  Kim, Oden and Song \cite{kikuchi-song-81, oden-kikuchi-82, oden-kim-82} (see also \cite{kikuchi-oden-88}) and the more recent study dealing with the boundary element method \cite{chernov-maischak-stephan-07}.

In the context of elliptic partial differential equations, the penalty method is also classical for the treatment of Dirichlet boundary conditions, and has been thoroughly analyzed for instance in \cite{babuska-73b,barrett-elliott-86}. However, due to the very different nature of Dirichlet boundary conditions and contact conditions, and of the resulting penalty methods, the aforementioned analysis can be hardly adapted to contact problems.

In this paper, we present a  convergence analysis of the penalty method for unilateral contact which use in particular some recent results from \cite{hild-renard-12}.
 We analyze both the continuous and discrete problems.  We limit the analysis to a conformal discretization with piecewise linear finite elements. In particular, we show that the same (sub-optimal, quasi-optimal or optimal) convergence rates as for the constrained problem (governed by a variational inequality) can be recovered with the choice $\varepsilon = h$, where $\varepsilon$ is the penalty parameter, and $h$ is the mesh size. A remarkable fact is that this choice is independent of the regularity of the continuous solution.

Let us introduce first some useful notations. In what follows, bold letters like $\u,\vv$, indicate vector or tensor valued quantities, while the capital ones (e.g., $\V,\K, \ldots$) represent functional sets involving vector fields. As usual, we denote by $(H^{s}(\cdot))^d$, $s\in \R, d=1,2,3$ the Sobolev spaces in one, two or three space dimensions (see \cite{adams-75}). The Sobolev norm of $(H^{s}(D))^d$ (dual norm if $s < 0$) is denoted by $\|\cdot\|_{s,D}$ and we keep the same notation when $d=1$, $d=2$ or $d=3$.  The letter $C$ stands for a generic constant, independent of the discretization parameters.

\section{Setting}
\label{sect:setting}

\subsection{The contact problem}
\label{sub:contact}

We consider an elastic body $\Omega$ in $\R^d$ with $d=2$ or $d=3$. Small
strain assumptions are made, as well as plane strain when $d=2$.
The boundary $\partial\Omega$ of
$\Omega$ is polygonal or polyhedral and we suppose that $\partial\Omega$
consists in three nonoverlapping parts $\Gamma_D$, $\Gamma_N$ and
the contact boundary $\Gamma_C$, with meas$(\Gamma_D) > 0$ and
meas$(\Gamma_C) > 0$. The contact boundary is supposed to be a
straight line segment when $d=2$ or a polygon when $d=3$ to simplify.
 The normal unit outward vector on
$\partial\Omega$ is denoted $\n$. In its initial stage, the
body is in contact on $\Gamma_C$ with a rigid foundation
and we suppose that the unknown
final contact zone after deformation will be included into
$\Gamma_C$. The body is clamped on $\Gamma_D$ for the sake of
simplicity. It is subjected to volume forces $\f \in
(L^2(\Omega))^d$ and to surface loads $\g \in
(L^2(\Gamma_N))^d$.

The unilateral contact problem in linear elasticity consists in
finding the displacement field $\u: \Omega \rightarrow \R^d$ verifying
the equations and conditions (\ref{eq})--(\ref{t.cont}):
\begin{align}
\label{eq}
\begin{aligned}
{\bf div}\, \bs(\u)  + \f &= \0   &\qquad \hbox{ in
} \Omega, \\
\bs(\u)&= {\bf A} \: \e (\u)  &\qquad
\hbox{ in } \Omega, \\
\u &= \0  &\qquad \hbox{ on
} \Gamma_{D}, \\
\bs(\u)  \n  &= \g &\qquad  \hbox{ on }
 \Gamma_N,
\end{aligned}
\end{align}
where $\bs = (\sigma_{ij}), \;1\le i,j \le d,$ stands for the
stress tensor field and ${\bf div}$ denotes the divergence
operator of tensor valued functions. The notation $\e(\vv) =
(\bnabla \vv + \bnabla \vv^{^T})/2$ represents the linearized strain
tensor field and $ {\bf A}$ is the fourth order symmetric
elasticity tensor having the usual uniform ellipticity and
boundedness property. For any displacement field $\vv$ and for any
density of surface forces $\bs(\vv) \n$ defined on $\partial
\Omega$ we adopt the following decomposition in normal and tangential components:
$$\vv = v_n \n +  \vv_t  \quad \hbox{ and } \quad
 \bs(\vv) \n = \sigma_{n}(\vv)\n + \bs_{t}(\vv).
$$
The nonlinear conditions describing unilateral contact  on
 $\Gamma_C$ are:
\begin{eqnarray}
 u_n  \le 0, \quad  \sigma_n(\u) \le 0,  \quad \sigma_n(\u)\, u_n
 =0,
\label{contact}
\end{eqnarray}
and the frictionless condition is
\begin{eqnarray}
\bs_t(\u)  = 0.
      \label{t.cont}
\end{eqnarray}
We introduce the Hilbert space:
$$
\V :=\left\{\vv\in \left(H^1(\Omega)\right)^d \; : \vv=\0 \hbox{ on }\Gamma_D\right\}.
$$
The convex cone of admissible displacements which satisfy the noninterpenetration
on the contact zone $\Gamma_C$ is:
$$
\K := \left\{\vv \in \V \::\: v_n = \vv \cdot \n \leq 0 \hbox{ on }\Gamma_C\right\}.
$$
Define
\begin{eqnarray*}
  a(\u, \vv) &:=& \int_{\Omega}  \bs(\u) :
  \e(\vv)
  \;d\Omega,\\
  L(\vv) &:=&  \int_{\Omega} \f \cdot \vv\;d\Omega +\int_{\Gamma_N} \g
  \cdot \vv\; \dg,
\end{eqnarray*}
for any $\u$ and $\vv$ in $\V$. From the previous assumptions, we know that $a(\cdot,\cdot)$
is bilinear, symmetric, $\V$-elliptic and continuous on $\V \times \V$. We know also that $L(\cdot)$
is a continuous linear form on $\V$.

The weak formulation of Problem \eqref{eq}-\eqref{t.cont}, as a variational inequality (see \cite{fichera-63} and also \cite{haslinger-hlavacek-80,haslinger-96,kikuchi-oden-88}), is:
\begin{equation} \label{weak}
\left\{
\begin{array}{l}
\hbox{Find } \u \in \K \hbox{ such that:}\\
a(\u,\vv - \u) \geq L(\vv -\u), \qquad \forall \, \vv \in \K.
\end{array}
\right.
\end{equation}
Stampacchia's Theorem ensures that Problem \eqref{weak} admits a unique solution.

\subsection{The penalty formulation of the unilateral contact problem}
\label{sub:penalty}

Let us first introduce the notation $[\cdot]_+$ for the positive part of a scalar quantity
$a \in \R$:
\[
[a]_+  =
\left \{ \begin{array}{cc}
a & \hbox{ if } a > 0,\\
0 & \hbox{ otherwise. }
\end{array}
\right.
\]
In the rest of this paper, we will make an extensive use of the following properties:
\begin{equation}
\label{plusprop}
a \leq [a]_+, \quad a[a]_+ = [a]^2_+, \quad \forall a \in \R.
\end{equation}
We introduce also:
$$W = \left\{ v_n \restrictiona{\Gc}: \vv \in \V \right\},$$
and its topological dual space $W'$, endowed with its usual dual norm.
Since $\Gc$ is a straight line segment or a polygon, we have
$H^{1/2}_{_{00}}(\Gc) \subset W \subset H^{1/2}(\Gc)$ which
implies $W' \subset (H^{1/2}_{_{00}}(\Gc))'$ where
$H^{1/2}_{_{00}}(\Gc)$ is the space of the
restrictions on $\Gc$ of functions in $H^{1/2}(\partial\Omega)$
vanishing outside $\Gc$ and
$H^{1/2}(\Gc)$
is the space of the restrictions on $\Gc$ of traces on $\partial
\Omega$ of functions in $H^1(\Omega)$.  We refer to \cite{lions-magenes-72} and
\cite{adams-75} for a detailed presentation of trace operators
and/or trace spaces.

Let $\eps > 0$ be a small parameter. The penalty method for
unilateral contact problem \eqref{weak} reads:

\begin{equation} \label{penalty}
\left\{
\begin{aligned}
&\hbox{Find } \u_\eps \in \V \hbox{ such that:}\\
&a(\u_\eps,\vv) + \frac1\eps \langle [u_{\eps,n}]_+ , v_n \rangle_{\Gamma_C} =
L(\vv), \qquad \forall \, \vv \in \V,
\end{aligned}
\right.
\end{equation}

where $\langle \cdot , \cdot \rangle_{\Gamma_C}$ stands for the duality product between
$W'$ and $W$.
Note that this formulation (which is an approximation of the exact contact conditions \eqref{contact}) is obtained by setting the condition $\sigma_n(\u_\eps) = - \frac1\eps [u_{\eps,n}]_+$ instead of the conditions in (\ref{contact}) on the contact boundary $\Gamma_C$.

\br
Note that for Problem \eqref{penalty}, it holds : $u_{\eps,n},v_n \in W \subset L^2(\Gamma_C)$, and
so $[u_{\eps,n}]_+ \in L^2(\Gamma_C)$. It results that the duality product $\langle \cdot , \cdot \rangle_{\Gamma_C}$ can also be understood simply as the inner product in $L^2(\Gamma_C)$. This remark still holds for the remaining part of the paper
 in which all the duality products $\langle \cdot , \cdot \rangle_{\Gamma_C}$ can be changed with inner products.
\er

We recall that Problem \eqref{penalty} is well-posed using an argument proposed by H. Brezis for
M type  and pseudo-monotone operators \cite{brezis-68} (see also \cite{lions-69} and \cite{kikuchi-song-81}):

\bt
\label{wellpenalty}
For all $\eps > 0$, Problem \eqref{penalty} admits a unique solution $\u_\eps$.
\et

\proof
Using the Riesz representation theorem, we define a (nonlinear) operator $\BB : \V \rightarrow \V$ with the following formula:

\[
( \BB \vv , \w )_{1,\Omega} := a (\vv , \w ) + \frac1\eps
\langle [v_{n}]_+ , w_n \rangle_{\Gamma_C}, \quad \forall \vv,\w \in \V,
\]
where $(.,.)_{1,\Omega} $ denotes the inner product in $(H^1(\Omega))^d$.
Note that Problem \eqref{penalty} is well-posed if and only if $\BB$ is a one-to-one operator.

Let $\vv,\w \in \V$, it follows from the definition of $\BB$ that:
\begin{align*}
(\BB \vv - \BB \w , \vv - \w)_{1,\Omega}
= a ( \vv - \w, \vv - \w) + \frac1\eps
\langle [v_{n}]_+ - [w_{n}]_+ , v_n - w_n \rangle_{\Gamma_C}.
\end{align*}
Due to the properties
\eqref{plusprop}, we observe that, for all $a,b\in\R$:
\begin{align*}
\begin{aligned}
([a]_+ - [b]_+)(a-b) & =  a[a]_+ + b[b]_+ -b[a]_+ -a[b]_+\\
& \geq  [a]^2_+ + [b]^2_+ - 2[a]_+[b]_+ \\ & =  ([a]_+ - [b]_+)^2 \geq 0.
\end{aligned}
%\label{plusprop2}
\end{align*}
This property combined to the $\V$-ellipticity of $a(\cdot,\cdot)$
imply that there exists $\alpha >0$ such that:
\begin{align}
\label{propB1}
(\BB \vv - \BB \w , \vv - \w)_{1,\Omega}
\geq \alpha \| \vv - \w \|^2_{1,\Omega}, \quad
\forall \vv,\w \in \V.
\end{align}
Let us also show that the operator $\BB$ is hemicontinuous, which means that for all
$\vv,\w \in \V$, the real function
\[
[0,1] \ni t \mapsto \varphi (t) := ( \BB ( \vv - t \w ) , \w )_{1,\Omega} \in \R
\]
is continuous. For $s,t \in [0,1]$, we have:
 \begin{align*}
| \varphi (t) - \varphi (s) | &
= | ( \BB ( \vv - t \w ) - \BB ( \vv - s \w ) , \w )_{1,\Omega} |\\
& \leq | a(\vv - t \w,\w) - a(\vv - s \w,\w) |
+ \frac1\eps | \langle [v_{n} -t w_n]_+ , w_n \rangle_{\Gamma_C}
- \langle [v_{n} -s w_n]_+ , w_n \rangle_{\Gamma_C} |\\
& =| a((s- t) \w,\w) |
+  \frac1\eps \left| \int_{\Gamma_C} \left([v_{n} -t w_n]_+ - [v_{n} -s w_n]_+\right) w_n \; d\Gamma \right| \\
& \leq | t - s | a(\w,\w)
+  \frac1\eps \int_{\Gamma_C}| [v_{n} -t w_n]_+ - [v_{n} -s w_n]_+| \; |w_n |\; d\Gamma.
\end{align*}

Using then the inequality $ | [v_n]_+ -  [w_n]_+ | \leq | v_n - w_n |$, we obtain:
\begin{align*}
| \varphi (t) - \varphi (s) | &  \leq | t - s | a(\w,\w)
+  \frac1\eps \int_{\Gamma_C}|t-s|  \;  |w_n | \; |w_n |\; d\Gamma \\
& = | t - s | \left ( a(\w,\w) + \frac1\eps \| w_n \|^2_{0,\Gamma_C}
\right) .
\end{align*}
It follows that $\varphi$ is Lipschitz, so continuous. The operator $\BB$ is then hemicontinuous.
Since \eqref{propB1} also holds,
we can apply the Corollary 15 (p.126) of \cite{brezis-68} to conclude that $\BB$ is a one-to-one
operator from $\V$ to $\V$. That concludes the proof of the theorem.
 \cqfd
\subsection{Finite element setting and discrete penalty problem}

Let $\V^h \subset \V$ be a family of finite dimensional vector
 spaces (see \cite{ciarlet-91}) indexed by $h$ coming from a family ${\cal T}^h$
of triangulations of the domain $\Omega$ ($h =
\max_{T \in {\cal T}^h} h_T$ where $h_T$ is the diameter of $T$). The
family of triangulations is supposed regular, i.e., there exists
$\sigma>0$ such that $\forall T \in {\cal T}^h,  h_T / \rho_T \le \sigma$
where $\rho_T$ denotes the radius of the inscribed ball %circle
in $T$.
We choose standard continuous and piecewise affine functions,
i.e.:
\begin{equation} \nonumber
\V^h = \left\{\vv^h \in (C(\overline\Omega))^d:
\vv^h\restrictiona{T} \in (\PP_1(T))^d, \forall T \in {\cal T}^h,
\vv^h = \0
\hbox{ on } \Gamma_D \right\}.
\end{equation}

We introduce $W^h(\Gamma_C)$, the space of normal traces on $\Gamma_C$ for discrete functions in $\V^h$:

\[
W^h(\Gamma_C) := \Big\{ \mu_h \in \caC(\overline{\Gamma_C}) \::\: \exists \vv^h \in \V^h, \vv^h \cdot \n = \mu_h \Big\}.
\]

To simplify, we suppose that the end points (or the border for $d=3$) of $\Gamma_C$ belong to $\Gamma_N$. Moreover we assume that the mesh on $\Gamma_C$ induced by ${\cal T}^h$, on which are defined functions of $W^h(\Gamma_C)$, is quasi-uniform, which implies in particular that it is
locally quasi-uniform, in the sense of Bramble \& al \cite{bramble-02}.

The discrete version of the penalty method \eqref{penalty} for Problem \eqref{weak} reads:

\begin{equation} \label{penaltymef}
\left\{
\begin{aligned}
&\hbox{Find } \u^h_{\eps} \in \V^h \hbox{ such that:}\\
&a(\u^h_\eps,\vv^h) + \frac1\eps \langle [u^h_{\eps,n}]_+ , v^h_n \rangle_{\Gamma_C} =
L(\vv^h), \qquad \forall \, \vv^h \in \V^h,
\end{aligned}
\right.
\end{equation}

where $\eps > 0$ is still the small penalty parameter.
Using exactly the same argument as in Theorem \ref{wellpenalty}, we see that
this problem admits one unique solution.

\section{Convergence analysis of the penalty method}
\label{sect:ca_analysis}

We present in this section the convergence analysis of both continuous and discrete penalty methods. We first state the two main theorems. Next, we give some necessary technical lemmas, followed by the proof of each theorem.

\subsection{Main results}

\subsubsection{Convergence when $\varepsilon \to 0$}

The penalty formulation introduces a consistency error. The following theorem gives a bound for this error as a function of the penalty parameter $\eps$. In particular, we recover the well-known result that the solution $\u_\eps$ of Problem \eqref{penalty} converges to the solution $\u$ of Problem
\eqref{weak} as $\eps \rightarrow 0$.
%Note that we do not recover the same result as in \cite{kikuchi-oden-88} since the latter reference involves a bound in %p 54, line 4  which we not use.
We suppose that $\u$ is more than $H^{\frac{3}{2}} $ regular  so that all the duality pairings become $L^2$ inner products.

\bt
\label{cvcont}
Suppose that $\u$, the solution of Problem \eqref{weak}, belongs to $(H^{\frac32+\nu} (\Omega))^d$ with $\nu \in (0,1/2]$. Let $\u_\eps$ be the solution of Problem \eqref{penalty}. We have the {\it a-priori} estimates:

\begin{align}
\| \u - \u_\eps \|_{1,\Omega}
&\leq C \eps^{\frac12+\nu} \| \u \|_{\frac32+\nu,\Omega}, \label{estim1}\\
 \left\| \sigma_n(\u) + \frac1\eps
[ u_{\eps,n} ]_+ \right\|_{0,\Gamma_C}
&\leq C \eps^{ \nu} \| \u \|_{\frac32+\nu,\Omega}, \label{estim2}\\
\left\| \sigma_n(\u) + \frac1\eps
[ u_{\eps,n} ]_+ \right\|_{-\nu,\Gamma_C}
&\leq C \eps^{2\nu} \| \u \|_{\frac32+\nu,\Omega} \label{estim3},
\end{align}

with $C > 0$ a constant, independent of $\eps$ and $\u$.

\et

\subsubsection{Convergence when $\varepsilon \to 0, h \to 0$}

Next, the following theorem provides the convergence rates for the discrete penalty problem, as a function of both discretization parameters $\eps$ and $h$.

\bt
\label{cvmef}
Suppose that the solution $\u$ of Problem \eqref{weak} belongs to
$(H^{\frac32+\nu} (\Omega))^d$ with $\nu \in (0, 1/2]$. The solution $\u^h_\eps$ of the discrete penalty problem
\eqref{penaltymef} satisfies the following error estimates in two space dimensions:

\begin{align}
\label{cvmef1}
\begin{aligned}
 \| \u - \u^h_\eps \|_{1,\Omega} & +  \eps^{\frac12}
 \left\| \sigma_n (\u) + \frac1\eps [u^h_{\eps,n}]_+ \right\|_{0,\Gamma_C} \\  &\leq C
 \left \{
\begin{array}{lll}
 \left ( h^{\frac12+\frac\nu2+\nu^2}  + h^{\nu} \eps^\frac12 + h^{\nu - \frac12} \eps
 \right )   \| \u \|_{\frac32+\nu,\Omega}& \quad \textrm{ if } 0<\nu<\frac12,\\ \\
 \left ( h | \ln h |^{\frac12} +  (h\eps)^{\frac12} + \eps \right )
 \| \u \|_{2,\Omega} & \quad \textrm{ if } \nu=\frac12,
\end{array}
\right.
\end{aligned}
\end{align}

\begin{align}
\label{cvmef2}
\begin{aligned}
& \left\| \sigma_n (\u) + \frac1\eps [u^h_{\eps,n}]_+ \right\|_{-\nu,\Gamma_C}
\\
 \leq  &\:\: C
 \left \{
\begin{array}{lll}
 \left (
 h^{\frac12+\frac{3\nu}{2}+\nu^2}  \eps^{-\frac12} +   h^{\frac{3\nu}{2}+\nu^2}
+ h^{2\nu - \frac12} \eps^{\frac12}  +  h^{2\nu-1} \eps
 \right )  \| \u \|_{\frac32+\nu,\Omega}& \textrm{ if } 0<\nu<\frac12,\\ \\
 \left (
    h^{\frac32}  \vert \ln h \vert^{\frac12}  \eps^{-\frac12} + h  \vert \ln h \vert^{\frac12}
 +  (h \eps)^{\frac12} % h^{\frac12} \eps^{\frac12}
 +    \eps
 \right )
 \| \u \|_{2,\Omega} & \textrm{ if } \nu=\frac12,
\end{array}
\right.
\end{aligned}
\end{align}
with $C > 0$ a constant, independent of $\eps$, $h$ and $\u$.
In three space dimensions, the terms $h^{\frac12+\frac\nu2+\nu^2}$ (resp. $h | \ln h |^{\frac12} $) in (\ref{cvmef1}) and (\ref{cvmef2})  have to be replaced with $h^{\frac12+\frac\nu2}$ (resp. $h^{\frac34} $).
\et

One interesting fact is that  the choice
$\eps=h$ in the previous theorem, leads to the same error estimates as for the finite element approximation of the variational inequality
 (see \cite{hild-renard-12}). With this choice, the error estimates are straightforward and are given next.

\bc
\label{cvopt}
 Suppose that the solution $\u$ of Problem \eqref{weak} belongs to
$(H^{\frac32+\nu} (\Omega))^d$ with  $\nu \in (0, 1/2]$.
Suppose also that the penalty parameter is chosen as $\eps = h$.
The solution $\u^h_\eps$ of the discrete penalty problem
\eqref{penaltymef} satisfies the following error estimates
in two space dimensions:

\begin{align}
\label{cvmef3}
\begin{aligned}
& \| \u - \u^h_\eps \|_{1,\Omega}
+ h^{\frac12}
\left\| \sigma_n (\u) + \frac1\eps [u^h_{\eps,n}]_+ \right\|_{0,\Gamma_C}
+  h^{{\frac12}-\nu}
\left\| \sigma_n (\u) + \frac1\eps [u^h_{\eps,n}]_+ \right\|_{-\nu,\Gamma_C} \\
& \leq  C
 \left \{
\begin{array}{ll}
  h^{\frac12+\frac\nu2+\nu^2}
  \| \u \|_{\frac32+\nu,\Omega}& \quad \textrm{ if } 0<\nu<\frac12,\\ \\
  h | \ln h |^{\frac12}
 \| \u \|_{2,\Omega} & \quad \textrm{ if } \nu=\frac12,
\end{array}
\right.
\end{aligned}
\end{align}

with $C > 0$ a constant, independent of  $h$ and $\u$. In three space dimensions, the terms $h^{\frac12+\frac\nu2+\nu^2}$ (resp. $h | \ln h |^{\frac12} $) in (\ref{cvmef3})   have to be replaced with $h^{\frac12+\frac\nu2}$ (resp. $h^{\frac34} $).
\ec

\br
Note that the  quasi-optimality is only due to the estimation of the contact term, and is also
present for standard finite element discretization \cite{hild-renard-12}. This is not an intrinsic
property of the penalty method. Note that optimality could be recovered in both 2D and 3D cases (and for any $0 < \nu \le 1/2)$ if additional assumptions on the transition area between contact and non contact are made (see \cite{hueber-wohlmuth-05}).
\er

\br
If $\u$  belongs to
$(H^{s} (\Omega))^2$ with  $s > 2$ ($\Omega \subset \R^2$) then it is easy to show (by using the result in \cite{renard-12}) that the error term in (\ref{cvmef3}) is bounded by $C h
 \| \u \|_{s,\Omega} $.
\er

\br To ensure convergence of the term $\| \u - \u^h_\eps \|_{1,\Omega} +  \eps^{\frac12}
 \left\| \sigma_n (\u) + \frac1\eps [u^h_{\eps,n}]_+ \right\|_{0,\Gamma_C}$, we need to choose at least $\eps < C h^{\frac12 - \nu}$ (see (\ref{cvmef1})). As stated previously, to recover optimal convergence rate of the terms (\ref{cvmef1}) and (\ref{cvmef2}), the best choice is $\eps = h$.
\er

\subsection{Preliminary technical lemmas}

Let $\caP^h : L^2(\Gamma_C) \rightarrow W^h(\Gamma_C)$ denote the $L^2(\Gamma_C)$-projection operator
onto $W^h(\Gamma_C)$. We recall in this lemma the stability and interpolation properties of
$\caP^h$:

\bl
Suppose that the mesh associated to $W^h(\Gamma_C)$ is locally quasi-uniform.
For all $s \in [0,1]$ and all $v\in H^s(\Gamma_C)$, we have the stability estimate:

\begin{equation}
\label{stabP}
\| \caP^h v \|_{s,\Gamma_C} \leq C \| v \|_{s,\Gamma_C}.
\end{equation}
The following interpolation estimate also holds:
\begin{equation}
\label{intP}
\| v - \caP^h v \|_{0,\Gamma_C} \leq C h^s \| v \|_{s,\Gamma_C},
\end{equation}
for all $v\in H^s(\Gamma_C)$. The constant $C > 0$ is in both cases independent of $v$ and $h$.
\el

\proof The stability estimate \eqref{stabP} is proven in \cite{bramble-02}. The interpolation estimate comes from e.g. \cite{bernardi-maday-patera-94}. \cqfd

%We recall a lemma from e.g. \cite{bernardi-maday-patera-94} which ensures the
We need another lemma which concerns the existence of a discrete bounded lifting from $\Gamma_C$ to $\Omega$:

\bl
Suppose that the mesh on the contact boundary $\Gamma_C$ is quasi-uniform.
There exists $\caR^h : W^h (\Gamma_C) \rightarrow \V^h$ and $C>0$, such that:

\begin{equation}
\label{lifting}
\caR^h(v^h)|_{\Gamma_C} \cdot \n = v^h, \quad
%\caR^h(v_1,v_2)|_{\Gamma_C} \cdot \t = v_2, \quad
\| \caR^h (v^h) \|_{1,\Omega} \leq C \| v^h \|_{\frac12,\Gamma_C},
\end{equation}

for all $v^h \in W^h (\Gamma_C)$.
\el

\proof The existence of such an operator is proven in \cite{bjorstad-widlund-86} (see also \cite{dominguez-sayas-03}).\cqfd

%\proof See e.g. \cite{bernardi-maday-patera-94}. \cqfd

The proof of each theorem relies strongly on appropriate estimates for the approximate contact condition on $\Gamma_C$ in dual norm. We give these estimates in the following lemma:

\bl
Suppose that
$\u$, the solution of Problem \eqref{weak}, belongs to $(H^{\frac32+\nu} (\Omega))^d$ with $\nu \in (0,1/2]$.
Let $\u_\eps$ be the solution of Problem \eqref{penalty}. We have the bound:

\begin{align}
\label{negboundcont}
\begin{aligned}
\left\| \sigma_n (\u) + \frac1\eps [u_{\eps,n}]_+ \right\|_{-\nu,\Gamma_C}
\leq  \: C \left (
\eps^\nu
\left\| \sigma_n (\u) + \frac1\eps [u_{\eps,n}]_+ \right\|_{0,\Gamma_C}
+
\eps^{\nu - \frac12} \:\| \u - \u_\eps \|_{1,\Omega}
 \right ),
\end{aligned}
\end{align}

with $C > 0$ a constant, independent of $\eps$, $\u$ and $\u_\eps$.

Let $\u^h_\eps$ be the solution of Problem \eqref{penaltymef}. We have also the bound:

\begin{align}
\label{negboundmef}
\begin{aligned}
\left\| \sigma_n (\u) + \frac1\eps [u^h_{\eps,n}]_+ \right\|_{-\nu,\Gamma_C}
\leq  \: C \left (
\: h^\nu \left\| \sigma_n (\u) + \frac1\eps [u^h_{\eps,n}]_+ \right\|_{0,\Gamma_C}
+
h^{\nu - \frac12} \:\| \u - \u^h_\eps \|_{1,\Omega}
 \right),
\end{aligned}
\end{align}

with $C > 0$ a constant, independent of $\eps$, $h$, $\u$ and $\u^h_\eps$.
\el

\proof
We first prove the bound \eqref{negboundmef} for the discrete solution.
The bound  \eqref{negboundcont} for the continuous solution will be derived in a quite similar manner.
By definition of the norm $\|\cdot\|_{-\nu,\Gamma_C}$ and using the projection operator $\caP^h$:

\begin{align}
&\left\| \sigma_n (\u) + \frac1\eps [u^h_{\eps,n}]_+ \right\|_{-\nu,\Gamma_C}  \nonumber \\  = &
\sup_{v \in H^{\nu} (\Gamma_C)}
\frac{\langle \sigma_n (\u) + \frac1\eps [u^h_{\eps,n}]_+ , v \rangle_{\Gamma_C}}{\| v \|_{\nu,\Gamma_C}}
  \nonumber  \\
 \le &
\sup_{v \in H^{\nu} (\Gamma_C)}
\frac{\langle \sigma_n (\u) + \frac1\eps [u^h_{\eps,n}]_+ , v - \caP^h v \rangle_{\Gamma_C}}{\| v \|_{\nu,\Gamma_C}} +
\sup_{v \in H^{\nu} (\Gamma_C)}
\frac{\langle \sigma_n (\u) + \frac1\eps [u^h_{\eps,n}]_+ , \caP^h v \rangle_{\Gamma_C}}{\| v \|_{\nu,\Gamma_C}}  \nonumber  \\
\leq &
\: \left\| \sigma_n (\u) + \frac1\eps [u^h_{\eps,n}]_+ \right\|_{0,\Gamma_C} \;
\sup_{v \in H^{\nu} (\Gamma_C)}
\frac{ \| v - \caP^h v \|_{0,\Gamma_C}}{\| v \|_{\nu,\Gamma_C}} + C
\sup_{v \in H^{\nu} (\Gamma_C)}
\frac{\langle \sigma_n (\u) + \frac1\eps [u^h_{\eps,n}]_+ , \caP^h v \rangle_{\Gamma_C}}{\|  \caP_h v \|_{\nu,\Gamma_C}} \nonumber \\
\leq & \: C \left (
\: h^\nu \left\| \sigma_n (\u) + \frac1\eps [u^h_{\eps,n}]_+ \right\|_{0,\Gamma_C}
+
 \sup_{v \in H^{\nu} (\Gamma_C)}
\frac{\langle \sigma_n (\u) + \frac1\eps [u^h_{\eps,n}]_+ , \caP_h v \rangle_{\Gamma_C}}{\| \caP_h v \|_{\nu,\Gamma_C}} \right).  \label{aa1}
\end{align}

In the fourth line, we used the Cauchy-Schwarz inequality and the stability property \eqref{stabP}. The last line is a direct consequence of the interpolation property \eqref{intP}.
Now, with help of the relationship $a(\u - \u^h_\eps , \vv^h) =
\langle \sigma_n (\u) + \frac1\eps [u^h_{\eps,n}]_+ , v^h_n \rangle_{\Gamma_C}$ for all $\vv^h
\in \V^h$, of the discrete lifting \eqref{lifting} and using the continuity of $a(\cdot,\cdot)$,
it results that:

\begin{align*}
\sup_{v \in H^{\nu} (\Gamma_C)}
\frac{\langle \sigma_n (\u) + \frac1\eps [u^h_{\eps,n}]_+ , \caP^h v \rangle_{\Gamma_C}}{\| \caP^h v \|_{\nu,\Gamma_C}}
= & \sup_{v \in H^{\nu} (\Gamma_C)}
\frac{\langle \sigma_n (\u) + \frac1\eps [u^h_{\eps,n}]_+ ,
\caR^h (\caP^h v)|_{\Gamma_C} \cdot \n  \rangle_{\Gamma_C}}{\| \caP^h v \|_{\nu,\Gamma_C}}\\
= & \sup_{v \in H^{\nu} (\Gamma_C)}
\frac{a(\u - \u^h_\eps , \caR^h (\caP^h v))}{\| \caP^h v \|_{\nu,\Gamma_C}}\\
\leq & \: C\| \u - \u^h_\eps \|_{1,\Omega}
\sup_{v \in H^{\nu} (\Gamma_C)}
\frac{\| \caR^h (\caP^h v) \|_{1,\Omega}}{\| \caP^h v \|_{\nu,\Gamma_C}} \\
\leq & \: C\| \u - \u^h_\eps \|_{1,\Omega}
\sup_{v \in H^{\nu} (\Gamma_C)}
\frac{\| \caP^h v \|_{\frac12,\Gamma_C}}{\| \caP^h v \|_{\nu,\Gamma_C}}.
\end{align*}

Since the discrete trace space $W^h(\Gamma_C)$ has the quasi-uniform mesh property,
we make use of the inverse inequality
\[
\| \caP^h v \|_{\frac12,\Gamma_C} \leq C h^{\nu - \frac12}
\| \caP^h v \|_{\nu,\Gamma_C},
\]

so that we finally get:

\begin{align*}
\sup_{v \in H^{\nu} (\Gamma_C)}
\frac{\langle \sigma_n (\u) + \frac1\eps [u^h_{\eps,n}]_+ , \caP^h v \rangle_{\Gamma_C}}{\| \caP^h v \|_{\nu,\Gamma_C}}
\leq \: C h^{\nu - \frac12} \:\| \u - \u^h_\eps \|_{1,\Omega}.
\end{align*}

This, together with (\ref{aa1}) leads to:

\begin{align*}
%\label{errestimneg}
\begin{aligned}
\left\| \sigma_n (\u) + \frac1\eps [u^h_{\eps,n}]_+ \right\|_{-\nu,\Gamma_C}
\leq  C \left (
\: h^\nu \left\| \sigma_n (\u) + \frac1\eps [u^h_{\eps,n}]_+ \right\|_{0,\Gamma_C}
+
h^{\nu - \frac12} \:\| \u - \u^h_\eps \|_{1,\Omega}
 \right),
\end{aligned}
\end{align*}

which is the desired bound (\ref{negboundmef}).

For the continuous problem, we introduce  $\V^\eps$, a fictitious finite element space, defined identically as $\V^h$ and with the choice of mesh size $h=\eps$.
We define also a fictitious discrete trace space $W^\eps (\Gamma_C)$, in the same manner
as $W^h(\Gamma_C)$, and with also the mesh size $h = \eps$.
We note simply $\caP^\eps : L^2(\Gamma_C) \rightarrow W^\eps (\Gamma_C)$ the $L^2(\Gamma_C)$-projection operator on $W^\eps (\Gamma_C)$. We again obtain an analogous estimate as in  (\ref{aa1}) :

\begin{align*}
%\label{normlessnu}
\begin{aligned}
& \left\| \sigma_n (\u) + \frac1\eps [u_{\eps,n}]_+ \right\|_{-\nu,\Gamma_C}
\\
\leq & \: C \left (
\eps^\nu \left\| \sigma_n (\u) + \frac1\eps [u_{\eps,n}]_+ \right\|_{0,\Gamma_C}
+
\sup_{v \in H^{\nu} (\Gamma_C)}
\frac{\langle \sigma_n (\u) + \frac1\eps [u_{\eps,n}]_+ , \caP^\eps v \rangle_{\Gamma_C}}{\| \caP^\eps v \|_{\nu,\Gamma_C}} \right).
\end{aligned}
\end{align*}

With help of the relationship
\begin{equation} \label{aa2}
a(\u - \u_\eps , \vv) =
\langle \sigma_n (\u) + \frac1\eps [\u_{\eps,n}]_+ , v_n \rangle_{\Gamma_C}, \qquad  \vv \in \V,
\end{equation}
using the continuity of $a(\cdot,\cdot)$ and a discrete lifting operator $\caR^\eps$, we come to the conclusion that:
\begin{align*}
%\label{normlessnu3}
\begin{aligned}
\left\| \sigma_n (\u) + \frac1\eps [u_{\eps,n}]_+ \right\|_{-\nu,\Gamma_C}
\leq  C \left (
\eps^\nu
\left\| \sigma_n (\u) + \frac1\eps [u_{\eps,n}]_+ \right\|_{0,\Gamma_C}
+
\eps^{\nu - \frac12} \:\| \u - \u_\eps \|_{1,\Omega}
 \right )
\end{aligned}
\end{align*}
which is the desired result (\ref{negboundcont}).\cqfd
\begin{remark}
When $\nu = 1/2$  the previous result can be improved, by using simply estimate (\ref{aa2}) to obtain
$$\left\| \sigma_n (\u) + \frac1\eps [u_{\eps,n}]_+ \right\|_{-1/2,\Gamma_C}
\leq  C \| \u - \u_\eps \|_{1,\Omega}.
 $$
\end{remark}

\subsection{Proof of Theorem \ref{cvcont}}

%\proof
We first use the $\V$-ellipticity of $a(\cdot,\cdot)$ and the fact that $\u$ (resp. $\u_\eps$)
is the solution of Problem \eqref{weak} (resp. \eqref{penalty}) to obtain:

\begin{align*}
\alpha \| \u - \u_\eps \|^2_{1,\Omega} &
\leq a ( \u - \u_\eps , \u - \u_\eps )\\
& = a ( \u , \u - \u_\eps ) - a ( \u_\eps , \u - \u_\eps )\\
& = \langle \sigma_n (\u) + \frac1\eps [u_{\eps,n}]_+ , u_n - u_{\eps,n} \rangle_{\Gamma_C}\\
&
= \langle \sigma_n (\u) , u_n  \rangle_{\Gamma_C}
+ \langle \frac1\eps [u_{\eps,n}]_+ , u_n \rangle_{\Gamma_C}
- \langle \sigma_n (\u) , u_{\eps,n} \rangle_{\Gamma_C}
- \langle \frac1\eps [u_{\eps,n}]_+ , u_{\eps,n} \rangle_{\Gamma_C},
\end{align*}

where $\alpha > 0$ denotes the ellipticity constant.
Contact conditions \eqref{contact} on $\Gamma_C$ yield:

\[
\langle \sigma_n (\u) , u_n \rangle_{\Gamma_C} = 0, \quad
\langle \frac1\eps [ u_{\eps,n} ]_+ , u_n \rangle_{\Gamma_C} \leq 0.
\]

Once again conditions \eqref{contact} associated to properties \eqref{plusprop} provide:

\[
- \langle \sigma_n(\u) ,
 u_{\eps,n} \rangle_{\Gamma_C} \leq
- \langle \sigma_n(\u) ,
 [u_{\eps,n}]_+  \rangle_{\Gamma_C},
\]

\[
\langle \frac1\eps [ u_{\eps,n} ]_+ ,
u_{\eps,n} \rangle_{\Gamma_C}
= \langle \frac1\eps [ u_{\eps,n} ]_+ ,
[ u_{\eps,n} ]_+ \rangle_{\Gamma_C}.
\]

We then get:

\begin{align*}
\alpha \| \u - \u_\eps \|^2_{1,\Omega} &
\leq
- \langle \sigma_n (\u) + \frac1\eps [u_{\eps,n}]_+ , [u_{\eps,n}]_+ \rangle_{\Gamma_C}.
\end{align*}

So we continue bounding:

\begin{align*}
%\label{pcv-calcint1}
\begin{aligned}
&\alpha \| \u - \u_\eps \|^2_{1,\Omega} \\
\leq & - \langle \eps (\sigma_n (\u) + \frac1\eps [u_{\eps,n}]_+) ,
\frac1\eps [u_{\eps,n}]_+ + \sigma_n (\u) - \sigma_n (\u) \rangle_{\Gamma_C}\\
= &
- \eps \left\| \sigma_n (\u) + \frac1\eps [u_{\eps,n}]_+ \right\|^2_{0,\Gamma_C}
+ \eps \langle \sigma_n (\u) + \frac1\eps [u_{\eps,n}]_+ , \sigma_n (\u) \rangle_{\Gamma_C}\\
\leq &
- \eps \left\| \sigma_n (\u) + \frac1\eps [u_{\eps,n}]_+ \right\|^2_{0,\Gamma_C}
+ \eps^{\delta} \left\| \sigma_n (\u) + \frac1\eps [u_{\eps,n}]_+ \right\|_{-\nu,\Gamma_C}
\: \eps^{1-\delta} \| \sigma_n (\u) \|_{\nu,\Gamma_C}\\
\leq &
- \eps \left\| \sigma_n (\u) + \frac1\eps [u_{\eps,n}]_+ \right\|^2_{0,\Gamma_C}
+ \frac{\eps^{2\delta}}{2\beta} \left\| \sigma_n (\u) + \frac1\eps [u_{\eps,n}]_+ \right\|^2_{-\nu,\Gamma_C}
+ \frac{\beta \eps^{2-2\delta}}2 \| \sigma_n (\u) \|^2_{\nu,\Gamma_C},
\end{aligned}
\end{align*}

with $\delta \in [0,1]$, $\beta > 0$. Note that since we supposed $\u \in (H^{\frac32+\nu}(\Omega))^2$,  we have
$\sigma_n (\u) \in H^{\nu} (\Gamma_C)$.
% We then derive a bound for
%$\| \sigma_n (\u) + \frac1\eps [\u_{\eps,n}]_+ \|_{-\nu,\Gamma_C}$ as follows:
We combine this result with estimation \eqref{negboundcont}:

\begin{align*}
%\label{pcv-calcint2}
\begin{aligned}
\alpha \| \u - \u_\eps \|^2_{1,\Omega}
\leq &
- \eps \left ( 1 - C \frac{\eps^{2(\delta+\nu) - 1}}{\beta} \right)
 \left\| \sigma_n (\u) + \frac1\eps [u_{\eps,n}]_+ \right\|^2_{0,\Gamma_C}\\
&+ C \frac{\eps^{2(\delta+\nu) - 1}}{\beta} \| \u - \u_\eps \|^2_{1,\Omega}
 + \frac{\beta \eps^{2-2\delta}}2 \| \sigma_n (\u) \|^2_{\nu,\Gamma_C},
\end{aligned}
\end{align*}

which can be transformed into:

\begin{align*}
%\label{pcv-calcint3}
\begin{aligned}
&\left( \alpha - C \frac{\eps^{2(\delta+\nu) - 1}}{\beta} \right)  \| \u - \u_\eps \|^2_{1,\Omega}
+ \eps \left ( 1 - C \frac{\eps^{2(\delta+\nu) - 1}}{\beta} \right)
 \left\| \sigma_n (\u) + \frac1\eps [u_{\eps,n}]_+ \right\|^2_{0,\Gamma_C}\\
\leq &
\quad \frac{\beta \eps^{2-2\delta}}2 \| \sigma_n (\u) \|^2_{\nu,\Gamma_C}.
\end{aligned}
\end{align*}

Taking $\delta = 1/2 - \nu$, $\beta = 2C\max(1,\alpha^{-1})$ and using the estimate
$\| \sigma_n (\u) \|_{\nu,\Gamma_C} \le C \| \u \|_{\frac32+\nu,\Omega}$
proves the bounds (\ref{estim1}) and (\ref{estim2}) of the theorem.
The bound (\ref{estim3}) of the theorem is then a direct consequence of this last result and
the estimation in norm of the negative exponent Sobolev space \eqref{negboundcont}.
\cqfd

\subsection{Proof of Theorem \ref{cvmef}}

We denote by $\caI^h$ the Lagrange interpolation operator mapping onto $\V^h$. We first use the
$\V$-ellipticity and the continuity of $a(\cdot,\cdot)$, as well as Young's inequality, to obtain:

\begin{align*}
\alpha \| \u - \u^h_\eps \|^2_{1,\Omega} &
\leq a ( \u - \u^h_\eps , \u - \u^h_\eps )\\
& =  a ( \u - \u^h_\eps , (\u - \caI^h \u) + (\caI^h \u - \u^h_\eps) )\\
& \leq C \| \u - \u^h_\eps \|_{1,\Omega} \| \u - \caI^h \u \|_{1,\Omega}
+ a ( \u - \u^h_\eps , \caI^h \u - \u^h_\eps )\\
& \leq \frac\alpha2 \| \u - \u^h_\eps \|^2_{1,\Omega}
+ \frac{C^2}{2\alpha} \| \u - \caI^h \u \|^2_{1,\Omega}
+ a ( \u  , \caI^h \u - \u^h_\eps )
- a ( \u^h_\eps , \caI^h \u - \u^h_\eps ),
\end{align*}

with $\alpha > 0$ the ellipticity constant. Since $\u$ is solution of \eqref{weak} and
$\u^h_\eps$  is solution of \eqref{penaltymef}, we can transform the term
$a ( \u  , \caI^h \u - \u^h_\eps ) - a ( \u^h_\eps , \caI^h \u - \u^h_\eps )$. So we obtain:

\begin{align*}
\frac\alpha2 \| \u - \u^h_\eps \|^2_{1,\Omega} &
\leq \frac{C^2}{2\alpha} \| \u - \caI^h \u \|^2_{1,\Omega} +
\langle \sigma_n(\u) + \frac1\eps [ u^h_{\eps,n} ]_+ ,
(\caI^h \u)_n - u^h_{\eps,n} \rangle_{\Gamma_C}.
\end{align*}

Because of conditions \eqref{t.cont} and since the 1D-Lagrange interpolation with piecewise-linear polynomials preserves the positivity,  we have that $(\caI^h \u)_n \leq 0$  on $\Gamma_C$. This implies:

\[
\langle \frac1\eps [ u^h_{\eps,n} ]_+ , (\caI^h \u)_n \rangle_{\Gamma_C} \leq 0.
\]

Once again condition \eqref{t.cont} associated to properties \eqref{plusprop}
yield:

\[
- \langle \sigma_n(\u) ,
 u^h_{\eps,n} \rangle_{\Gamma_C} \leq
- \langle \sigma_n(\u) ,
 [u^h_{\eps,n}]_+  \rangle_{\Gamma_C},
\]

\[
\langle \frac1\eps [ u^h_{\eps,n} ]_+ ,
u^h_{\eps,n} \rangle_{\Gamma_C}
= \langle \frac1\eps [ u^h_{\eps,n} ]_+ ,
[ u^h_{\eps,n} ]_+ \rangle_{\Gamma_C}.
\]

This results into:

\begin{align}
\label{errestim1}
\frac\alpha2 \| \u - \u^h_\eps \|^2_{1,\Omega} &
\leq \frac{C^2}{2\alpha} \| \u - \caI^h \u \|^2_{1,\Omega} +
\langle \sigma_n(\u) , (\caI^h \u)_n \rangle_{\Gamma_C} -
\langle \sigma_n(\u) + \frac1\eps [ u^h_{\eps,n} ]_+ ,
 [u^h_{\eps,n}]_+ \rangle_{\Gamma_C}.
\end{align}

We bound the last term
$\langle \sigma_n(\u) + \frac1\eps [ u^h_{\eps,n} ]_+ ,
 [u^h_{\eps,n}]_+ \rangle_{\Gamma_C}$ as follows:

\begin{align*}
& - \langle \sigma_n(\u) + \frac1\eps [ u^h_{\eps,n} ]_+ ,
 [u^h_{\eps,n}]_+ \rangle_{\Gamma_C} \\ = & -
\langle \eps (\sigma_n (\u) + \frac1\eps [u^h_{\eps,n}]_+) ,
\frac1\eps [u^h_{\eps,n}]_+ + \sigma_n (\u) - \sigma_n (\u) \rangle_{\Gamma_C}\\
 = &
- \eps \left\| \sigma_n (\u) + \frac1\eps [u^h_{\eps,n}]_+ \right\|^2_{0,\Gamma_C}
+ \eps \langle \sigma_n (\u) + \frac1\eps [u^h_{\eps,n}]_+ , \sigma_n (\u) \rangle_{\Gamma_C}\\
\leq &
- \eps \left\| \sigma_n (\u) + \frac1\eps [u^h_{\eps,n}]_+ \right\|^2_{0,\Gamma_C}
+ \eps \left\| \sigma_n (\u) + \frac1\eps [u^h_{\eps,n}]_+ \right\|_{-\nu,\Gamma_C}
\| \sigma_n (\u) \|_{\nu,\Gamma_C}.
\end{align*}

The last inequality has been made possible due to the assumption on the regularity of $\u$, which implies that
$\sigma_n (\u) \in H^{\nu} (\Gamma_C)$. Treating the last term with Young's inequality,
and with $\beta > 0$, we obtain finally:

\begin{equation}
\label{errestim2}
\begin{aligned}
& - \langle \sigma_n(\u) + \frac1\eps [ u^h_{\eps,n} ]_+ ,
 [u^h_{\eps,n}]_+ \rangle_{\Gamma_C}\\
& \leq
- \eps \left\| \sigma_n (\u) + \frac1\eps [u^h_{\eps,n}]_+ \right\|^2_{0,\Gamma_C}
+ \frac{\eps^{2}}{2\beta} \left\| \sigma_n (\u) + \frac1\eps [u^h_{\eps,n}]_+ \right\|^2_{-\nu,\Gamma_C}
+ \frac{\beta}2 \| \sigma_n (\u) \|^2_{\nu,\Gamma_C}.
\end{aligned}
\end{equation}

We now combine this last inequality \eqref{errestim2} with \eqref{errestim1}
and then insert the estimation \eqref{negboundmef} into the resulting inequality:

\begin{align*}
& \frac\alpha2 \| \u - \u^h_\eps \|^2_{1,\Omega}\\
\leq &\: \frac{C^2}{2\alpha} \| \u - \caI^h \u \|^2_{1,\Omega} +
\langle \sigma_n(\u) , (\caI^h \u)_n \rangle_{\Gamma_C} \\
&
- \eps \left\| \sigma_n (\u) + \frac1\eps [u^h_{\eps,n}]_+ \right\|^2_{0,\Gamma_C}
+ \frac{\eps^{2}}{2\beta} \left\| \sigma_n (\u) + \frac1\eps [u^h_{\eps,n}]_+ \right\|^2_{-\nu,\Gamma_C}
+ \frac{\beta}2 \| \sigma_n (\u) \|^2_{\nu,\Gamma_C}\\
\leq & \: \frac{C^2}{2\alpha} \| \u - \caI^h \u \|^2_{1,\Omega} +
\langle \sigma_n(\u) , (\caI^h \u)_n \rangle_{\Gamma_C} \\
&
- \eps \left ( 1 - C \frac{\eps h^{2\nu}}{2\beta} \right )
\left\| \sigma_n (\u) + \frac1\eps [u^h_{\eps,n}]_+ \right\|^2_{0,\Gamma_C}
+ C h^{2\nu - 1} \frac{\eps^{2}}{2\beta} \| \u - \u_\eps \|^2_{1,\Omega} + \frac{\beta}2 \| \sigma_n (\u) \|^2_{\nu,\Gamma_C}.
\end{align*}

We rearrange the terms:

\begin{align}
\label{errestim3}
\begin{aligned}
& \left ( \frac\alpha2 - C h^{2\nu - 1} \frac{\eps^{2}}{2\beta} \right ) \| \u - \u^h_\eps \|^2_{1,\Omega}
+ \eps \left ( 1 - C \frac{\eps h^{2\nu}}{2\beta}  \right )
\left\| \sigma_n (\u) + \frac1\eps [u^h_{\eps,n}]_+ \right\|^2_{0,\Gamma_C}\\
 \leq  & \:\: \frac{C^2}{2\alpha} \| \u - \caI^h \u \|^2_{1,\Omega} +
\langle \sigma_n(\u) , (\caI^h \u)_n \rangle_{\Gamma_C}
+ \frac{\beta}2 \| \sigma_n (\u) \|^2_{\nu,\Gamma_C}.
\end{aligned}
\end{align}

We choose $\beta$ as follows:

\[
\beta = C \left ( h^{2\nu - 1} \eps^2 \alpha^{-1}+ h^{2\nu} \eps \right ),
\]

with $C > 0$ a constant sufficiently large so that the two left terms in \eqref{errestim3}
are positive
irrespectively of the values of $\eps$ and $h$.
It results that:

\begin{align}
\label{estimbeta}
&\beta \| \sigma_n (\u) \|^2_{\nu,\Gamma_C} =  C \left ( h^{2\nu - 1} \eps^2 + h^{2\nu} \eps \right ) \| \sigma_n (\u) \|^2_{\nu,\Gamma_C}.
\end{align}

The estimation of the Lagrange interpolation error in $L^2$ and $H^1$ norms on a domain $D$ is classical (see e.g., \cite{ern-guermond-04}):

\begin{equation}
h^{-1}\| \u - \caI^h \u \|_{0,D}  +  \| \u - \caI^h \u \|_{1,D} \leq C h^{s-1}  \| \u \|_{s,D},  \label{estimint}
\end{equation}
for $s \in (1,2]$.
The contact term $\langle \sigma_n(\u) , (\caI^h \u)_n \rangle_{\Gamma_C}$
can be estimated in two space dimensions using results from \cite{hild-renard-12}:

\begin{equation}
\langle \sigma_n(\u) , (\caI^h \u)_n \rangle_{\Gamma_C} \leq C \label{2D}
\left \{
\begin{array} {lll}
h^{1+\nu+2\nu^2} \| \u \|^2_{\frac32+\nu,\Omega} & \quad \textrm{ if } 0<\nu<\frac12, \\
h^2 | \ln h | \| \u \|^2_{2,\Omega} & \quad \textrm{ if } \nu=\frac12.
\end{array}
\right.
\end{equation}

In three space dimensions the bound is obtained in a straightforward way using (\ref{estimint}) for any $ 0<\nu \le 1/2$ :

\begin{equation}
\langle \sigma_n(\u) , (\caI^h \u)_n \rangle_{\Gamma_C} \leq C
h^{1+\nu} \| \u \|^2_{\frac32+\nu,\Omega}.  \label{3D}
\end{equation}

In two space dimensions, we combine finally %these last three estimations with \eqref{errestim3}
the estimations \eqref{errestim3}$-$\eqref{2D}
to prove that:

\begin{align*}
%\label{errestim4}
\begin{aligned}
&  \| \u - \u^h_\eps \|^2_{1,\Omega}
+ \eps
\left\| \sigma_n (\u) + \frac1\eps [u^h_{\eps,n}]_+ \right\|^2_{0,\Gamma_C}\\
 \leq  &\:\: C
 \left \{
\begin{array}{lll}
 \left ( h^{1+\nu+2\nu^2} + h^{2\nu - 1} \eps^{2} + h^{2\nu} \eps
 \right )   \| \u \|^2_{\frac32+\nu,\Omega}& \quad \textrm{ if } 0<\nu<\frac12,\\ \\
 \left ( h^2 | \ln h | + \eps^{2} + h \eps  \right )
 \| \u \|^2_{2,\Omega} & \quad \textrm{ if } \nu=\frac12,
\end{array}
\right.
\end{aligned}
\end{align*}

which is the first required estimate.
This, together with the estimate \eqref{negboundmef}, yields additionally the bound in two space dimensions:

\begin{align*}
%\label{errestim5}
\begin{aligned}
& \left\| \sigma_n (\u) + \frac1\eps [u^h_{\eps,n}]_+ \right\|_{-\nu,\Gamma_C}
\\
 \leq  &\:\: C
 \left \{
\begin{array}{lll}
 \left ( (h^\nu \eps^{-\frac12} + h^{\nu-\frac12} ) h^{\frac12+\frac\nu2+\nu^2}
+ h^{2\nu - \frac12} \eps^{\frac12} +  h^{2\nu} +  h^{2\nu-1} \eps
 \right )  \| \u \|_{\frac32+\nu,\Omega}& \textrm{ if } 0<\nu<\frac12,\\ \\
 \left ( (h^{\frac12} \eps^{-\frac12} + 1 ) h \vert \ln h \vert^\frac12
 +  h^{\frac12} \eps^{\frac12} +  h +  \eps  \right )
 \| \u \|_{2,\Omega} & \textrm{ if } \nu=\frac12.
\end{array}
\right.
\end{aligned}
\end{align*}

Using  $h^{2\nu} \le h^{\frac{3\nu}2 + \nu^2}$ and $h \le  h \vert \ln h \vert^\frac12$ ends the proof of (\ref{cvmef2}). The bounds in three space dimensions are obtained as before using estimate (\ref{3D}) instead of  (\ref{2D}). \cqfd

%%%%%%%%%%%%%%%%%%%%%%%%%%%%%%%%%%%%%%%%%%%%%%%%
%{\bf Acknowledgements.}

%% The Appendices part is started with the command \appendix;
%% appendix sections are then done as normal sections
%% \appendix

%% \section{}
%% \label{}

%% References
%%
%% Following citation commands can be used in the body text:
%% Usage of \cite is as follows:
%%   \cite{key}          ==>>  [#]
%%   \cite[chap. 2]{key} ==>>  [#, chap. 2]
%%   \citet{key}         ==>>  Author [#]

%% References with bibTeX database:

\bibliographystyle{model1b-num-names}
\bibliography{biblio}

%% Authors are advised to submit their bibtex database files. They are
%% requested to list a bibtex style file in the manuscript if they do
%% not want to use model1-num-names.bst.

%% References without bibTeX database:

% \begin{thebibliography}{00}

%% \bibitem must have the following form:
%%   \bibitem{key}...
%%

% \bibitem{}

% \end{thebibliography}

\end{document}